\documentclass[preprint,12pt]{elsarticle}

\usepackage{amssymb}
\usepackage{amsmath}

\usepackage[T1]{fontenc} 
\usepackage{tikz}
\usepackage{hyperref}

\usepackage{amsthm}
\theoremstyle{plain}
\newtheorem{theorem}{Theorem}
\newtheorem{lemma}[theorem]{Lemma}
\newtheorem{proposition}[theorem]{Proposition}

\theoremstyle{definition}

\newtheorem{conjecture}[theorem]{Conjecture}
\newtheorem{observation}[theorem]{Observation}
\numberwithin{equation}{section}

\newcommand{\N}{\mathbb N}
\newcommand{\A}{\mathcal A}
\newcommand{\B}{\mathcal B}

\def\uu{\mathbf{u}}
\def\xx{\mathbf{x}}
\def\vv{\mathbf{v}}
\def\CR{\mathrm{E}}

\journal{European Journal of Combinatorics}

\begin{document}

\begin{frontmatter}



\title{On two conjectures of Shallit\\ about Thue-Morse-like sequences}


\author[1]{Lubom\'{i}ra Dvo\v{r}\'{a}kov\'{a}\corref{cor}}
\ead{lubomira.dvorakova@fjfi.cvut.cz}

\author[2]{Savinien  Kreczman\fnref{skfunding}}
\ead{savinien.kreczman@uliege.be}

\author[1]{Edita Pelantov\'{a}}
\ead{edita.pelantova@fjfi.cvut.cz}

\cortext[cor]{Corresponding author.}
\fntext[skfunding]{Supported by the FNRS Research Fellow grant 1.A.789.23F, by a Fédération Wallonie-Bruxelles Travel Grant and by a WBI Excellence WORLD grant. Avec le soutien de Wallonie-Bruxelles International. \includegraphics[scale=0.01]{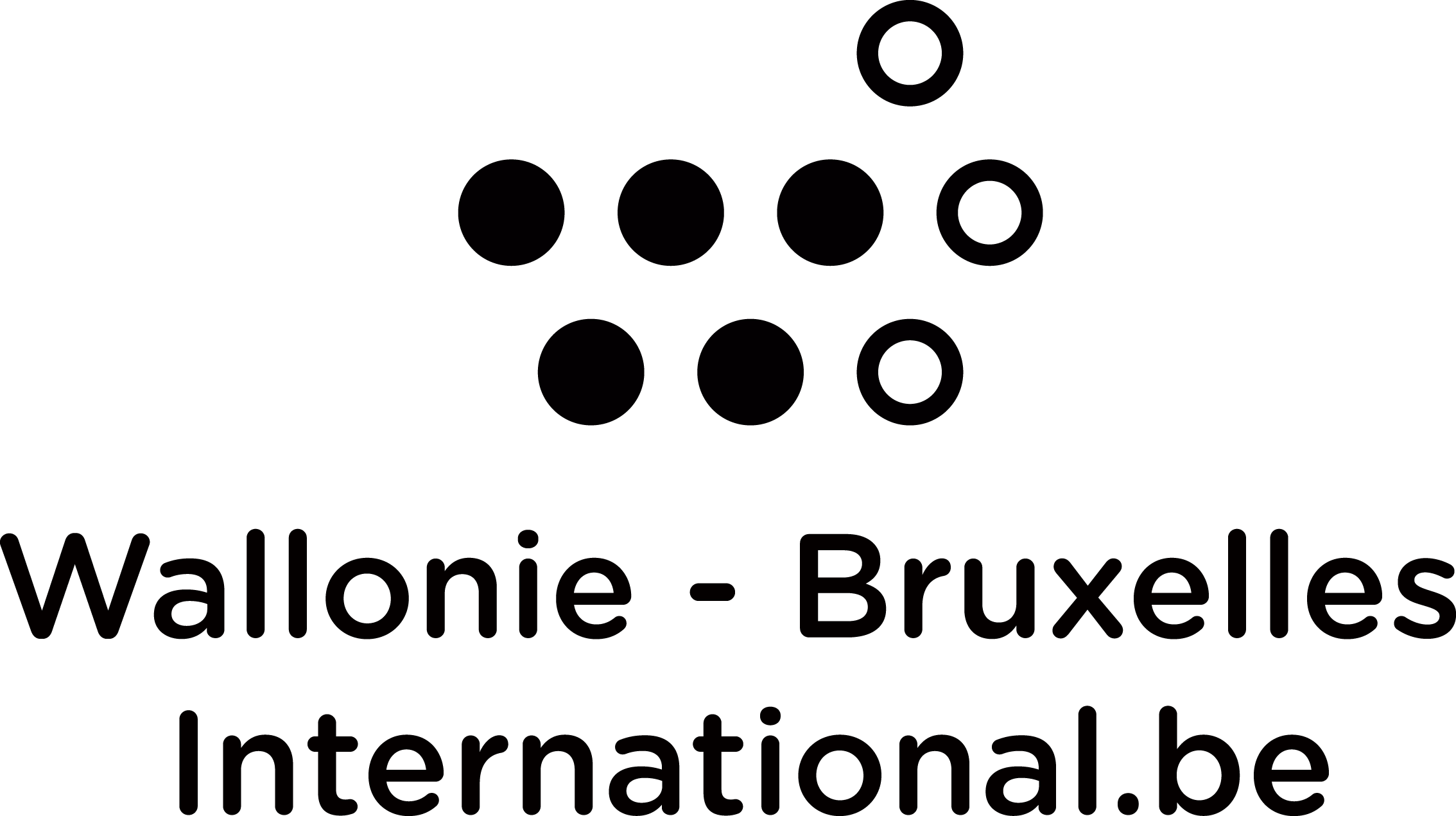}}

\affiliation[1]{organization={Department of Mathematics, Faculty of Nuclear Sciences and 
Physical Engineering, Czech Technical University in Prague},
            addressline={Trojanova 13}, 
            city={Praha 2},
            postcode={120 00}, 
            country={Czech Republic}}
\affiliation[2]{organization={Department of Mathematics, University of Li\`{e}ge},
            addressline={\\ All\'{e}e de la D\'{e}couverte 12}, 
            city={Li\`{e}ge},
            postcode={4000}, 
            country={Belgium}}

\begin{abstract}
We study a class of infinite words $\xx_k$, $k \in \N, k\geq 1$, recently introduced by J. Shallit. This class includes the Thue-Morse sequence $\xx_1$, the 
Fibonacci-Thue-Morse sequence $\xx_2$,
and the Allouche-Johnson sequence $\xx_3$.  Shallit stated and for $k=3$ proved two conjectures  on properties of $\xx_k$.  The first  conjecture concerns the factor complexity, the second one  the critical exponent of these words. We confirm the validity of both  conjectures for every $k$. 
\end{abstract}



\begin{keyword}
morphism \sep Thue-Morse sequence \sep factor complexity \sep critical exponent \sep asymptotic critical exponent \sep bispecial factors



\end{keyword}

\end{frontmatter}

%
%
%

\section{Introduction}

In a recent paper \cite{Sha25}, Shallit leverages \texttt{Walnut} to prove a number of properties related to the Narayana morphism $\mu\colon 0\mapsto 01,\, 1\mapsto 2,\, 2\mapsto 0$, its fixed point the Narayana word, and the Narayana numeration system. The latter is the positional numeration system with weights $(U_n)$ satisfying the recurrence relation $U_{n+3}=U_{n+2}+U_n$ with initial conditions $1,2,3$.

This numeration system is one instance of a family defined by the recurrence relation $U_{n+k}=U_{n+k-1}+U_n$ with the initial conditions $1,\ldots,k$. The case $k=1$ corresponds to the usual binary numeration system, and $k=2$ to the Zeckendorf numeration system~\cite{Zec72}. This family of numeration systems is linked to the family of recursively defined Hofstadter functions, see \cite{MeekVanRees84,Let25}.

Shallit \cite{Sha25} established a connection to a sequence introduced by Allouche and Johnson \cite{AlJo2020}: the $n$-th term of this sequence is the parity of the number of $1$'s in the Narayana representation of $n$. Again, this is the $k=3$ representative of a family which contains the Thue-Morse sequence itself ($k=1$) and the Fibonacci-Thue-Morse sequence (\cite{Sha21}, $k=2$). We call other representatives of this family \emph{Thue-Morse-like sequences}, although Shallit simply calls them $\xx_k$.

Shallit shows that 
for $k\geq 2$, the sequence ${\xx}_{k}$ is  the image under the projection $\pi$ of the fixed point $\uu_k=\xi_k^\omega(0)$ of the morphism  $\xi = \xi_k$ defined over the alphabet $\mathcal{A} = \{0,1,\ldots, k-1, 
0',1',\ldots, (k-1)'\}$   as follows 
\begin{equation}\label{eq:xi}
    \xi_k: \qquad 
    \begin{array}{rclcccrcl}
        0&\to& 01 &&&& 0'&\to &0'1'\\
        1&\to& 2 &&&& 1'&\to& 2'\\
        2&\to& 3 &&&& 2'&\to& 3'\\
        &\vdots& &&\text{and}&& &\vdots&\\
        (k-2)& \to & (k-1)&&& & (k-2)'& \to & (k-1)'\\
        (k-1)&\to&0' & &&&(k-1)'&\to &0.
    \end{array}
\end{equation}
The projection $\pi=\pi_k$ we apply to  $\uu_k$ is given  by 
\begin{equation}\label{eq:projection}
    \pi(a)=
    \begin{cases}
        0, & \text{if } a =0 \text{ or }  a\in \{1', 2', \ldots, (k-1)'\};\\
        1, & \text{if } a =0'  \text{ or } a\in \{1, 2, \ldots, k-1\}.
    \end{cases}
\end{equation}

Having proven results for the Allouche-Johnson sequence $\xx_3$ and comparing to the known results for $\xx_1$ and $\xx_2$, Shallit formulates two conjectures regarding the family $\xx_k$.
The first conjecture concerns factor complexity, the second one critical exponent. See Section \ref{sec:prelim} for definitions of these notions.

\begin{conjecture}[\cite{Sha25}, Conjecture 45]\label{conj:complexity}
    Let $k \geq 1$. The first difference of factor complexity of ${\xx}_k$, for $n$ large enough, takes the values $4k-2$ and $4k$ only.
\end{conjecture}

\begin{conjecture}[\cite{Sha25}, Conjecture 44]\label{conj:repetitions}
    Let $k \geq 1$. The sequence ${\xx}_k$ has critical exponent $k + 1$, which is
    attained by the words $0^{k+1}$ and ${1}^{k+1}$. It contains no factor of length $2n + k$ and period $n$,
    and therefore has asymptotic critical exponent equal to $2$.
\end{conjecture}

The famous Thue-Morse sequence $\xx_1$ is known to satisfy both conjectures. Already Thue showed that $\xx_1$ is overlap-free and contains squares of unbounded size, see~\cite{Berstel2007}. Its critical exponent equals 2, which is minimal among binary sequences. Its asymptotic critical exponent also equals 2 and the first difference of factor complexity takes only the values 2 and 4 for positive arguments, see~\cite{Brlek1989} and~\cite{LuVa1989}. 
The cases of $k=2$ (the Fibonacci-Thue-Morse sequence ${\xx}_2$) and $k=3$ (the Allouche-Johnson sequence ${\xx}_3$) were studied already (\cite{Sha21}, resp. \cite{Sha25}). Both conjectures hold for these particular cases. However, the methods of the latter paper cannot be used to prove these conjectures for all $k$.

The main goal of this article is to prove both conjectures for all $k$. For this purpose, we first study bispecial factors in the fixed point $\uu_k$ of $\xi_k$ in Section \ref{sec:fixedPoint} and we show that the language of $\uu_k$ is overlap-free.  Section \ref{sec:propertiesOfProjection} 
examines properties of the Thue-Morse-like sequence $\xx_k= \pi(\uu_k)$, again with a focus on bispecial factors, and prepares tools to prove both conjectures.
The proofs  are presented in Sections \ref{sec:proof_complexity} (Conjecture 1) and~\ref{sec:proof_repetitions} (Conjecture 2).

\section{Preliminaries}\label{sec:prelim}

An \textit{alphabet} $\mathcal A$ is a finite set, its elements are \textit{letters}. A \textit{word} $u$ over $\mathcal A$ of \textit{length} $n$ is a finite sequence $u = u_0 u_1 \cdots u_{n-1}$ of letters $u_j\in\mathcal A$ for all $j \in \{0,1,\dots, n-1\}$. The length of $u$ is denoted $|u|$.
The set of all finite words over $\A$ is denoted $\A^*$. The set $\A^*$ equipped with concatenation as the operation forms a monoid with the \textit{empty word} $\varepsilon$ as the neutral element. 
Consider $u, p, s, v \in \A^*$ such that $u=pvs$, then the word $p$ is called a \textit{prefix}, the word $s$ a \textit{suffix} and the word $v$ a~\textit{factor} of $u$. 
The \textit{longest common prefix} of two words $u$ and $z$ is denoted $\text{lcp}\{u,z\}$, it is the longest factor of $u$ and $z$ that is a prefix of both $u$ and $z$. The \textit{longest common suffix} ($\text{lcs}$) is defined analogously.
A~\textit{sequence} (or an \textit{infinite word}) $\vv$ over $\A$ is an infinite sequence $\vv = v_0 v_1 v_2 \cdots$ of letters $v_j \in \A$ for all $j \in \N$. A~\textit{word} $w$ over $\mathcal A$ is called a~\textit{factor} of the sequence $\vv = v_0 v_1 v_2 \cdots$ if there exists $j \in \mathbb N$ such that $w = v_j v_{j+1} v_{j+2} \cdots v_{j+|w|-1}$. The index $j$ is called an \textit{occurrence} of $w$. If $j=0$, then $w$ is a \textit{prefix} of $\vv$. Consider a factor $w$ of an infinite word $\vv = v_0 v_1 v_2 \cdots$. Let $j < \ell$ be two consecutive occurrences of $w$ in $\uu$. Then the word $v=v_j v_{j+1} \cdots v_{\ell-1}$ is a \textit{return word} to $w$ in $\vv$ and $vw$ is a~\textit{complete return word} to $w$ in $\vv$.

The \textit{language} $\mathcal{L}(\vv)$ of a sequence $\vv$ is the set of factors occurring in $\vv$.
The \textit{factor complexity} of a sequence $\vv$ is the map ${\mathcal C}:\mathbb N \to \mathbb N$, where 
$$
{\mathcal C}(n)=\#\{w\in {\mathcal L}(\vv)\,: |w| = n\}.
$$ 

To evaluate the factor complexity, Cassaigne~\cite{Ca1997} introduced the concept of special factors.  
A~factor $w$ of a sequence $\vv$ is \textit{left special} if $aw, bw \in \mathcal{L}(\vv)$ for at least two distinct letters ${a,b} \in \A$. Such letters $a,b$ are called \textit{left extensions} of $w$. The set of all left extensions of $w$ in $ \vv$ is denoted  ${\rm Lext}(w)$.
A \textit{right special} factor,  \textit{right extensions} and ${\rm Rext}(w)$ are defined analogously. The factor $w$ is \textit{bispecial} if it is both left special and right special.

The \textit{first difference of factor complexity} is a mapping $\Delta \mathcal{C}:\mathbb N \to \mathbb N$ defined by $\Delta \mathcal{C}(n)= \mathcal{C}(n+1)- \mathcal{C}(n)$.
The first difference of factor complexity of a binary sequence ${\bf v} $ satisfies
\begin{equation}\label{eq:ComplexityDifference}
    \Delta \mathcal{C}(n) = \#\{w \in \mathcal{L}({\bf v}):  w  \text{\ is left special in } {\bf v} \text{\ and\ } |w| =n\}\,.
\end{equation}

A sequence $\vv$ is \textit{recurrent} if each factor of $\vv$ has infinitely many occurrences in $\vv$. Moreover, a recurrent sequence $\vv$ is \textit{uniformly recurrent} if for every $n \in \mathbb N$ there exists $N\in \mathbb N$ such that each factor of $\vv$ of length $N$ contains all factors of $\vv$ of length $n$. 

A sequence $\vv$ is \textit{eventually periodic} if there exist words $u$ (possibly empty) and $v$ such that $\vv=uv^{\omega}$, where $\omega$ denotes an infinite repetition. The sequence $\vv$ is called \textit{aperiodic} otherwise.
Let  $\vv$ be a~uniformly recurrent and aperiodic sequence. Then  every left special factor $v$  in $\vv$ is a prefix of a bispecial factor in $\vv$. There exists a unique shortest bispecial factor having a prefix $v$. 

A \textit{morphism} is a map $\psi: \A^* \to \B^*$ such that $\psi(uv) = \psi(u)\psi(v)$ for all words $u, v \in \A^*$.
The morphism $\psi$ can be naturally extended to a sequence $\vv=v_0 v_1 v_2\cdots$ over $\A$ by setting
$\psi(\vv) = \psi(v_0) \psi(v_1) \psi(v_2) \cdots\,$.
If a morphism $\psi: \A^* \to \A^*$ satisfies $\psi(a) \not =\varepsilon$ for every letter $a \in \mathcal A$ and $\psi(b)=bw$ for some $b \in \A$ and $w\in\A^*, w\not =\varepsilon$, then there exists a sequence $\vv$ having the prefix $\psi^n(b)$ for every $n \in \mathbb N$, thus $\psi(\vv)=\vv$ and $\vv$ is a~\textit{fixed point} of $\psi$. 
In the sequel, we use the notation $\psi^{\omega}(b)$.
A~morphism $\psi: \A^* \to \A^*$ is \textit{primitive} if there exists $k\in \mathbb N$ such that $\psi^k(a)$ contains all letters of $\mathcal A$ for each letter $a \in \mathcal A$. It is known that a fixed point of a primitive morphism is uniformly recurrent~\cite{Qu1987}.

Let ${\bf u}$ be a sequence over an alphabet $\mathcal A$ and let $\psi:{\mathcal A}^* \to {\mathcal B}^*$ be a~morphism. 
Consider a factor $w$ of $\psi({\bf u})$. We say that $(w_1, w_2)$ is a \emph{synchronizing point} of $w$ if $w=w_1w_2$ and for all $p,s \in {\mathcal L}(\psi({\bf u}))$ and $v \in {\mathcal L}({\bf u})$ such that $\psi(v)=pws$ there exists a factorization $v=v_1v_2$ of $v$ with $\psi(v_1)=pw_1$ and $\psi(v_2)=w_2s$. We denote the synchronizing point by $w_1${\tiny{\textbullet}}$w_2$.

Let $u$ be a word of length $q=|u|\geq 1$.  If $z$ of length $p$ is a~prefix of $u^\omega$, we say that $q$ is a~\textit{period} of $z$  and write $z = u^r$, where $r={\frac{p}{q}}$.   For instance, the Czech word \texttt{ka\v cka}~\footnote{ka\v cka: Czech crown (informally), duck (in a Slovak dialect)} has period 3 and $\texttt{ka\v cka} = (\texttt{ka\v c} )^{\frac{5}{3}}$.  The \emph{critical exponent} of a sequence $\vv$ is defined as 
$$ 
\CR(\vv)  = \sup \{ r\in \mathbb{Q} :   u^r\ \text{is a non-empty factor of} \ \vv\}\,. 
$$
Its asymptotic version is called the \textit{asymptotic critical exponent}, denoted $\CR^*(\vv)$. If $\CR(\vv)=\infty$, then $\CR^*(\vv)=\infty$. Otherwise, it is defined as
$$ 
\CR^*(\vv)  = \lim_{n \to \infty} \sup \{ r \in \mathbb{Q} : u^r\ \text{is a non-empty factor  of period $\geq n $ of} \ \vv\}\,. 
$$

To compute the (asymptotic) critical exponent of a uniformly recurrent aperiodic sequence, it suffices to know bispecial factors and their shortest return words.
\begin{theorem}[\cite{DolceDP2023}, Theorem 3]\label{thm:FormulaForE}
    Let $\vv$ be a uniformly recurrent aperiodic sequence.
    Let $(w_n)_{n\in\N}$ be the sequence of all bispecial factors in $\vv$ ordered by length.
    For every $n \in \N$, let $r_n$ be the shortest return word to the bispecial factor $w_n$ in $\vv$.
    Then
    $$
    E(\vv) = 1 + \sup\left\{\frac{|w_n|}{|r_n|} \ : \ n \in \mathbb N\right\}\quad \text{and}\quad E^*(\vv) = 1 + \limsup_{n \to \infty}\frac{|w_n|}{|r_n|}\,.
    $$
\end{theorem}

\section{Bispecial factors of the fixed point $\uu=\xi^\omega(0)$} \label{sec:fixedPoint}
Let us fix $k\geq 4$ as Conjectures~\ref{conj:complexity} and~\ref{conj:repetitions} are proven for $k<4$. We denote $\uu={\uu}_k$ the fixed point of $\xi=\xi_k$, defined in~\eqref{eq:xi}. 
Our first goal is to study bispecial factors and their shortest return words. This knowledge will enable us in the sequel to determine the first difference of complexity of $\xx=\pi(\uu)$ (Conjecture~\ref{conj:complexity}) and using Theorem~\ref{thm:FormulaForE} we will show that there are no overlaps in $\uu$, which is essential for the proof of Conjecture~\ref{conj:repetitions}.

\subsection{Properties of the fixed point $\uu$}
We denote $\mathcal A$ the alphabet, i.e., ${\mathcal A}=\{0,1,2,\dots, (k-1), 0',1',2',\dots, (k-1)'\}$.
In our proofs we will often use symmetry of the language of $\uu$. 
\begin{observation}
    The language of $\uu$ is closed under the letter permutation $\ell \leftrightarrow \ell'$ for each $\ell \in \{0,1,\dots, k-1\}$. Applying  this exchange  to every letter occurring in $w\in \mathcal{L}(\uu)$, we  obtain the factor $w'\in \mathcal{L}(\uu)$ and call it  {\em twin} of $w$.
\end{observation}

Let us list some useful basic observations that follow from the form of the morphism $\xi$, defined in~\eqref{eq:xi}.

\begin{observation}\label{obs:basic} 
    The following statements apply also for twins.
    \begin{enumerate}
    \item\label{item:synchr_letters} The synchronizing points of letters are {\tiny{\textbullet}}$0, 1${\tiny{\textbullet}}, {\tiny{\textbullet}}$2${\tiny{\textbullet}}, $\dots$, {\tiny{\textbullet}}$(k-1)${\tiny{\textbullet}}.
    \item $\xi^j(0)=01\cdots j$ for $j \in \{1,2, \ldots, k-1\}$.
    \item $\xi^k(0)=01\cdots (k-1)0'=\xi^{k-1}(0)0'$.
    \item\label{item:mk} If $\xi^n(0)$ ends in $0$ or $0'$, then $n=mk$ for some $m\in \mathbb N$.
    \item\label{item:1} $\xi^{k-1}(1)=0'$, hence $\xi^n(1)$ starts in $0'$ for each $n \geq k-1$.
    \item\label{item:j} The letter $j$ occurs in $\uu$ only as a suffix of $01\cdots j$ for $j \in \{1,2, \ldots, k-1\}$. 
    
    \item\label{item:predek0} Except for its first occurrence, the letter $0$ occurs in $\uu$ only as a suffix of the following factors: 
    
    --\quad $012\cdots \ell 0$,\, \quad where $\ell\in  \{1,2, \ldots, (k-1)\}$;  
     
    --\quad $0'1'2'\cdots \ell'0$, \quad where $\ell'\in  \{1',2', \ldots, (k-1)'\}$; 
    
    --\quad $ 012\,\cdots (k-1)0'0$; 
     
    --\quad $0'1'2'\cdots (k-1)'00$.
    
    \item\label{item:001} The word $00$ occurs in $\uu$ only as a prefix of $001$.
    
    \item\label{item:(k-1)0} The word $(k-1)0$ occurs in $\uu$ only as a prefix of $(k-1)01$. 
    
    \item\label{item:Rext} Let $w \in \mathcal{L}(\uu)$, $w\neq \varepsilon$, and $\#{\rm Rext}(w)\geq  2$. Then $\#{\rm Rext}(w)\leq 3$ and $0$ or $0'$ belong to ${\rm Rext}(w)$. If $\#{\rm Rext}(w)=3$, then $ \{0, 0'\}\subset{\rm Rext}(w) $. 
    \end{enumerate}
\end{observation}

\subsection{Klouda's method for determining bispecial factors} 

To describe the structure of all bispecial factors, we use a method by Klouda \cite{Klouda2012}. 
This method does not work directly with bispecial factors, but with bispecial triplets. 

The first part of a triplet $T$  is an element of the so-called left forky set $F_L$, the last part of the triplet is an element of the so-called right forky set $F_R$ and the middle part is a bispecial factor itself.

\begin{description}
\item[Forky sets] \quad The left forky set $F_L$ is formed by a finite number of unordered pairs of words $\{s_L, t_L\}$ from $\mathcal{A}^*$. Analogously, the elements of the right forky set  $F_R$ are denoted $\{s_R, t_R\}$.    
Consult Definition 20 in \cite{Klouda2012} to see that for our substitution $\xi$ we can take the forky set  $F_R$ and $F_L$ as follows. 

\begin{itemize}

    \item  The right forky set  $F_R$ is formed by pairs $\{s_R, t_R\}$ of two kinds:   

The first kind of pairs together with their twins are as follows (we consider all $\ell  \in \{1,2,\ldots, k-2\}$):
$$
\begin{array}{llll} 
    \{0,(k\!-\!1)'0\},&  \{01,(k\!-\!1)'0'\},&   \{00',(k\!-\!1)'0'\},&  \{00,(k\!-\!1)'0'\},
\\  
    \{0, \ell\},& \{0, \ell'\,0'\},&\{0, \ell'\,0\},&  \{0, \ell'\, (\ell\!+\!1)'\},
\\
    \{0,0'\},&\{0, k-1\}\,;&&
\end{array}
$$

The second kind of pairs in $F_R$ together with their twins are as follows (we  consider all $\ell ,m  \in \{1,2,\ldots, k-1\}$):
$$
\begin{array}{l}
    \{\ell,m'\},\\
    \{\ell 0,m0\},\quad \{\ell(\ell +1),m\}, \quad  \{\ell 0,m0'\},\quad \{\ell 0',m\} \quad  \text{for $\ell <m$},\\
    \{\ell 0, m (m+1)\} \quad \text{for $\ell <m\leq k-2$}\,.
\end{array}
$$

    \item $F_L = \bigl\{\{a,b\} : a,b \in \mathcal{A}, a\neq b\bigr\}$. 

\end{itemize}
Let us stress two important properties of $F_R$: 
the  first letter of $s_R$ is distinct from the first letter of $t_R$ and  for every pair $\{s_R,t_R\}\in F_R$  there exists a unique pair $\{ \tilde{s}_R, \tilde{t}_R\} \in F_R$ such that $p\tilde{s}_R$ is a prefix of $\xi(s_R)$ and  $p\tilde{t}_R$ is a prefix of $\xi(t_R)$ with 
$p = {\rm lcp}\{\xi(s_R), \xi(t_R)\}$. The unique pair $\{ \tilde{s}_R, \tilde{t}_R\}$ is denoted $f_R\bigl(\{s_R,t_R\}\bigr)$.  Note that for every  pair $\{s_R, t_R\} \in F_R$ of the second kind,  ${\rm lcp}\{\xi(s_R),\xi(t_R)\}\bigr)= \varepsilon$.
    
The set $F_L$ satisfies the analogous properties which allow us to define $f_L\bigl(\{s_L,t_L\}\bigr)$.  Here,   ${\rm lcs}\{\xi(s_L), \xi(t_L)\}$ is empty for every pair from $F_L$.

\item[Bispecial triplets]\quad    Given the sets $F_L$ and $F_R$, we call \emph{bispecial triplet} a triplet $T=\bigl(\{s_L, t_L\}, w, \{s_R, t_R\}\bigr)$ such that
\begin{itemize}
    \item $\{s_L, t_L\}\in F_L$ and   $\{s_R, t_R\}\in F_R$;  
    \item 
both $s_L w\,s_R$ and $t_L w \,t_R$ are factors of $\uu$. 
Obviously, $w$ is a bispecial factor of $\uu$.
\end{itemize}

\item[Transformation of triplets] \quad  
We define a mapping  $\Phi$  on the set of bispecial triplets
$$
\Phi(\{s_L, t_L\}, w, \{s_R, t_R\})=(\{\tilde{s}_L, \tilde{t}_L\}, \tilde{w}, \{\tilde{s}_R, \tilde{t}_R\})\,,
$$
where
$
\{\tilde{s}_L, \tilde{t}_L\} = f_L\bigl(\{s_L,t_L\}) 
$, $
\{\tilde{s}_R, \tilde{t}_R\} = f_R\bigl(\{s_R,t_R\}) 
$  and 
\begin{equation}\label{eq:middleBS}
    \tilde{w} = {\rm lcs}\{\xi(s_L), \xi(t_L)\} \, \varphi(w)\, {\rm lcp}\{\xi(s_R), \xi(t_R)\}. 
\end{equation}
Clearly, $\Phi(\{s_L, t_L\}, w, \{s_R, t_R\})$ is a bispecial triplet, too. 
\end{description}

\begin{theorem}(\cite{Klouda2012}, Theorem~36) 
    Let $\uu$ be a fixed point of a primitive morphism $\varphi$ over $\mathcal{A}$ and let $F_L$ be a~left forky set and $F_R$ a~right forky set. Then there exists a~finite set $\mathcal{I}$ of bispecial triplets such that every bispecial factor of $\uu$ occurs as the middle element of a~triplet $\Phi^n(T)$ for some $n \in \N$ and $T\in \mathcal{I}$.
\end{theorem}

As follows from the proof of the previous theorem, 
to find the set of initial triplets $\mathcal{I}$, it suffices to consider only triplets where the bispecial factor $w$ does not have a synchronizing point. For our substitution $\xi$, the only bispecial factor without the synchronizing point is the empty word. Hence,  the set $\mathcal{I}$ contains only triplets of the form $T_I=\bigl(\{s_L, t_L\}, \varepsilon, \{s_R, t_R\}\bigl)$.  
In our specific case of substitution $\xi$, the bispecial factor associated with the triplet $\Phi(T_I)$ by \eqref{eq:middleBS} equals $\tilde{w}={\rm lcs}\{\xi(s_L), \xi(t_L)\} \, \xi(\varepsilon)\, {\rm lcp}\{\xi(s_R), \xi(t_R)\}= {\rm lcp}\{\xi(s_R), \xi(t_R)\}$.  Thus, only elements of $F_R$ for which ${\rm lcp}\{\xi(s_R), \xi(t_R)\}\neq \varepsilon$ give rise to non-empty bispecial factors. All such pairs or their twins are listed below.   
\begin{equation}\label{eq:TypyBS} 
    \{0, (k\!-\!1)'0\},\quad \{01, \, (k\!-\!1)'0'\},  \quad  \{00', (k\!-\!1)'0'\} \quad \text{and} \quad 
    \{00, \, (k\!-\!1)'0'\} \,.
\end{equation}

By successively applying the transformation $\Phi$ to the above initial triplets, we will see in the next subsection that only elements $\{s_R,t_R \} \in F_R$  of the first kind appear in $\Phi^n(T_I)$.

\subsection{List of bispecial triplets in $\uu$}\label{subsec:triplets}
Now we can find all bispecial factors of $\uu$. 
As we explained, in the pair $\{s_R,t_R\}$ of any initial triplet (up to the twin relation) the first digit of $s_R$ must be $0$ and the first digit of $t_R$ must be $(k\!-\!1)'$. Hence the last digit of $t_L$ must be $(k\!-\!2)'$. We will discuss four types of initial triplets $T$ with the middle element $w =\varepsilon$ according with  the pairs $\{s_R,t_R\}$ listed  in \eqref{eq:TypyBS}.  

We can observe for any initial bispecial triplet $T$ that after a finite number of steps the bispecial triplet $\Phi^n(T)$ has the element from the right forky set $\{s_R, t_R\}=\{0,0'\}$. In other words, every sufficiently long bispecial factor of $\uu$ has right extensions $0$ and $0'$. Consequently, if $w$ is the bispecial factor associated with the bispecial triplet $\Phi^n(T)$, then $\xi(w)$ is the bispecial factor associated with the bispecial triplet $\Phi^{n+1}(T)$.
  
In the role of $\{s_L,t_L\}$ we sometimes write $\{*, a\}$. By this notation we mean a pair where one element is  $a \in \mathcal{A}$ and the other element  is any letter of $ \mathcal{A}$ different from $a$.

\begin{description} 
\item[Type 1 ] \ \ starting  with  \ \ 
$T=\bigl(\{*,(k-2)'\}, \ \varepsilon, \  \{0,(k-1)'0\}\bigr) $:\quad  
$$
\begin{array}{rcl}
    \Phi({T})&=&\bigl(\{*,(k-1)'\}, \ 0, \  \{1,0\}\bigr)\\
    \Phi^2({ T})&=& \bigl(\{*,0\},\ \xi(0),\  \{2,0\}\bigr)\\
    &\vdots \\
    \Phi^{k-1}({ T})&=& \bigl(\{*, k-3\},\ \xi^{k-2}(0),\  \{k-1,0\}\bigr)\\
    \Phi^{k}(T) &=& \bigl(\{*, k-2\},\ \xi^{k-1}(0),\  \{0',0\}\bigr).\\
\end{array}
$$
          
\item[Type 2]  \ \ starting with \ \ 
${T}=\bigl(\{*, (k-2)'\}, \ \varepsilon, \  \{01, (k-1)'0'\}\bigr)$:
$$
\begin{array}{rcl}
    \Phi({T})&=&\bigl(\{*,(k-1)'\},\  0,\  \{12,0'\}\bigr)\\
    \Phi^2({T})&=& \bigl(\{*, 0\},\ \xi(0), \ \{23,0'\}\bigr) \\
    &\vdots \\
    \Phi^{k-1}({T})&=& \bigl(\{*, k-3\},\ \xi^{k-2}(0),\ \{(k-1)0',0'\}\bigl)\\
    \Phi^{k}({T}) &=& \bigl(\{*, k-2\},\ \xi^{k-1}(0)0'= \xi^k(0),  \ \{0',1'\}\bigr) \\
    &\vdots \\
    \Phi^{2k-2}({ T}) &=& \bigl(\{*, (k-4)'\},\ \xi^{2k-2}(0), \  \{0',(k-1)'\}\bigr)\\
    \Phi^{2k-1}({ T}) &=& \bigl(\{*, (k-3)'\},\ \xi^{2k-1}(0), \  \{0',0\}\bigr).\\
\end{array}
$$
           
\item[Type 3] \ \ starting with  \ \ $T=(\{(k-1)', (k-2)'\}, \varepsilon, \{00',(k-1)'0'\})$:
$$
\begin{array}{rcl}
    \Phi(T)&=&\bigl(\{0,(k-1)'\}, \  0,  \ \{10',0'\}\bigr)\\
    \Phi^2(T)&=&\bigl(\{1,0\},\  \xi(0), \  \{20',0'\}\bigr)\\ 
    &\vdots \\
    \Phi^{k-1}(T)&=&\bigl(\{k-2,k-3)\}, \ \xi^{k-2}(0),  \ \{ (k-1)0',0'\}\bigr)\\[2mm] 
    & & \text{coincides with $ \Phi^{k-1}(T)$ in Type 2}.
\end{array}
$$
                  
\item[Type 4]\ \ starting with \ \  $T=\bigl(\{(k-1)', (k-2)'\}, \ \varepsilon, \  \{00,(k-1)'0'\}\bigr)$:
$$
\begin{array}{rcl}
    \Phi(T)&=&\bigl(\{0,(k-1)'\}, \ 0, \  \{10,0'\}\bigr)\\
    \Phi^2(T)&=&\bigl(\{1, 0\}, \ \xi(0), \ \{20,0'\}\bigr)\\ 
    &\vdots \\
    \Phi^{k-2}(T)&=&\bigl(\{k-3, k-4\}, \  \xi^{k-3}(0), \  \{(k-2)0, 0'\}\bigr)\\ 
    \Phi^{k-1}(T)&=&\bigl(\{k-2,k-3\}, \ \xi^{k-2}(0), \ \{ (k-1)0, 0'1'\}\bigr)\\ 
    \Phi^{k}(T)&=&\bigl(\{k-1, k-2\}, \ \xi^{k-1}(0)0'=\xi^k(0), \  \{0, 1'2'\}\bigr)\\ 
    &\vdots \\
    \Phi^{2k-2}({ T})&=&\bigl(\{(k-3)',(k-4)'\},\ \xi^{2k-2}(0), \  \{0,(k-1)'0\}\bigr)\\
    \Phi^{2k-1}({ T})&=& \bigl(\{(k-2)', (k-3)'\},\ \xi^{2k-1}(0)0, \ \{1,0\}\bigr)\\
    &\vdots \\
    \Phi^{3k-2}({ T})&=& \bigl(\{k-3, k-4\},\ \xi^{3k-2}(0)\xi^{k-1}(0), \  \{0',0\}\bigr).\\
\end{array}
$$
\end{description}
   
\begin{proposition}\label{prop:aboutBSfactors}  
    Every bispecial factor  of $\uu$ or its twin is equal to 
    $$
    u=\xi^n(0)\quad \text{or}\quad v=\xi^{n+2k-1}(0)\xi^{n}(0) \ \text{\ \ \ for some $n \in \N$.}
    $$
    Moreover, $\mathrm{Lext}(u)=\mathcal A$, while $\#\mathrm{Lext}(v)=2$.
\end{proposition}

\begin{proof}
The statement follows from the fact that all bispecial factors are associated with bispecial triplets. According to the description of bispecial triplets, see Types 1 to 4 (do not forget to consider twins),  we know that $\xi^n(0)$ has all letters as left extensions and $\xi^{n+2k-1}(0)\xi^{n}(0)$, coming only from Type 4, has two left extensions. 
\end{proof}

\subsection{The sequence $\uu$ is overlap-free}
Let us first state a useful lemma on properties of the images of $0$.

\begin{lemma}\label{lem:images0}
    Consider the sequence $\uu=\xi^{\omega}(0)$.
    \begin{enumerate}
    \item\label{item:synchr}$\xi^n(0)$ has synchronizing points {\tiny{\textbullet}}\hspace{0.5mm}$\xi^n(0)${\tiny{\textbullet}} for each $n \in \mathbb N, n \geq 1, n\not =k$.
    \item\label{item:retword} The shortest return word to $\xi^n(0)$ is $\xi^n(0)$ for each $n \leq k$.
    \item\label{item:next_letter} $\xi^k(0)\xi^k(0)$ is always followed by $0'$.
    \end{enumerate}
\end{lemma}

\begin{proof}
\begin{enumerate}
\item If $\xi^n(0)$ ends in $a\not \in \{0,0'\}$, then the statement is clear by Item~\ref{item:synchr_letters} of Observation~\ref{obs:basic}. If $\xi^n(0)$ ends in $0$ or $0'$, then $n=mk$ for some $m\in \mathbb N$ by Item~\ref{item:mk} of Observation~\ref{obs:basic}. Observing the list of bispecial triplets, we can see that $\text{Rext}(\xi^n(0))=\{0,0'\}$ for $n\geq 2k-1$. This implies that for $m \geq 2$, the only right extensions of $\xi^{mk}(0)$ are $0$ and $0'$, but then there are synchronizing points  {\tiny{\textbullet}}\hspace{0.5mm}$\xi^{mk}(0)${\tiny{\textbullet}}$0$ or  {\tiny{\textbullet}}\hspace{0.5mm}$\xi^{mk}(0)${\tiny{\textbullet}}$0'$. 
\item $\xi^n(00)$ occurs in $\uu$, so $\xi^n(0)$ is a return word to itself. By Item 2 and 3 of Observation~\ref{obs:basic} the factor $\xi^n(0)$ has all distinct letters. If there is a return word to $\xi^n(0)$ shorter than $\xi^n(0)$, then $\xi^n(0)$ is a (possibly fractional) power of this word, but then $\xi^n(0)$ cannot consist of mutually distinct letters. Consequently, $\xi^n(0)$ is the shortest return word to $\xi^n(0)$.
\item By Observation~\ref{obs:basic} we have  $\xi^k(0)\xi^k(0)=\xi^{k-1}(0)0'\xi^{k-1}(0)0'$. By Item~\ref{item:synchr} $\xi^{k-1}(0)$ occurs only as the $(k-1)$-st image of $0$. Therefore $\xi^k(0)\xi^k(0)=\xi^{k-1}(010)0'$ occurs only as a prefix of the $(k-1)$-st image of $0101$ or $0100'$. However, $0100'$ is not a factor of $\uu$ and $0101$ has a unique right extension $2$ (by Item~\ref{item:001} of Observation~\ref{obs:basic}), hence  $\xi^k(0)\xi^k(0)$ occurs only as a prefix of $\xi^{k-1}(01012)=\xi^k(0)\xi^k(0)0'1'$.
\end{enumerate}
\end{proof}

To show that $\uu$ does not contain overlaps, we will apply Theorem~\ref{thm:FormulaForE}. For this purpose, we need some more knowledge on return words to bispecial factors. 

\begin{lemma}\label{lem:retwords}
    The shortest return word $r$ to a bispecial factor $w$ in $\uu$ satisfies:
    \begin{itemize}
        \item $r=w$ if $w=\xi^n(0)$ for $n \in \mathbb N$;
        \item $|r|\geq |w|$ if $w=\xi^{n+2k-1}(0)\xi^{n}(0)$ for $n\in \mathbb N$\,.
    \end{itemize}
\end{lemma}

\begin{proof}
For $w=\xi^n(0)$, where $0\leq n\leq k$, Item~\ref{item:retword} of Lemma~\ref{lem:images0} implies that the shortest return word to $w$ equals $w$ (moreover, each return word to $w$ contains $w$ as its prefix).
If $n>k$, then $w$ has synchronizing points {\tiny{\textbullet}}$w${\tiny{\textbullet}}\hspace{0.5mm} by Item~\ref{item:synchr} of Lemma~\ref{lem:images0}. 
If $rw$ is a complete return word to such bispecial factor $w$, then $rw=\xi(r'\xi^{n-1}(0))$, where $r'$ is a return word to $\xi^{n-1}(0)$. Hence $r=\xi(r')$.
This proves the first claim.

To get the second claim, let us explain that the factor $\xi^{2k-1}(0)$ is not a~return word to $\xi^{2k-1}(0)0$ because $\xi^{2k-1}(0)\xi^{2k-1}(0)$ is always followed by $0'$. By the same arguments as above, we get that $\xi^{2k-1}(0)\xi^{2k-1}(0)$ occurs only as a prefix of the $(k-1)$-st image of $\xi^k(0)\xi^k(0)$. By Item~\ref{item:next_letter} of Lemma~\ref{lem:images0}, $\xi^k(0)\xi^k(0)$ is always followed by $0'$, which implies that $\xi^{2k-1}(0)\xi^{2k-1}(0)$ is always followed by $0'$. Consequently, $\xi^{2k-1}(0)0$ is a prefix of any return word to $\xi^{2k-1}(0)0$. Thus, using Item~\ref{item:synchr} of Lemma~\ref{lem:images0}, every return word to $w=\xi^{n+2k-1}(0)\xi^n(0)$, where $n\not=k$, is the $n$-th image of $\xi^{2k-1}(0)0$ and has therefore $w$ as prefix. Moreover, for $n=k$, the factor $w=\xi^{3k-1}(0)\xi^k(0)=\xi^{3k-1}(0)\xi^{k-1}(0)0'$ occurs only as the prefix of the $(k-1)$-st image of $\xi^{2k}(0)00'$ or $\xi^{2k}(0)01$. Since the last letter of $\xi^{2k}(0)$ is $0$, the word $\xi^{2k}(0)00'$ is not a factor of $\uu$ by Item~\ref{item:001} of Observation~\ref{obs:basic}. Consequently, 
$w=\xi^{3k-1}(0)\xi^k(0)$ also occurs only as the $k$-th image of $\xi^{2k-1}(0)0$ and each return word to $w$ has $w$ as its prefix. The second claim then follows.
\end{proof}

\begin{proposition}\label{pro:NoOverlap}
    The sequence $\uu=\xi^{\omega}(0)$ does not contain overlaps and contains squares of unbounded size. Hence, the critical and asymptotic critical exponent of $\uu$ equal $2$.
\end{proposition}

\begin{proof}
There exist squares of unbounded size since $(\xi^n(0))^2$ is a factor of $\uu$ for any $n\in \mathbb N$.

To show that there are no overlaps, using Theorem~\ref{thm:FormulaForE}, it is sufficient to show that $\frac{|w|}{|r|}\leq 1$ for any bispecial factor $w$ of $\uu$ and its shortest return word $r$. This follows from Lemma~\ref{lem:retwords}.
\end{proof}

\section{Properties of the projection $\pi$}\label{sec:propertiesOfProjection}

Let us recall that we fix $k\geq 4$. We denote $\xi=\xi_k$ and $\pi=\pi_k$ the morphisms defined in~\eqref{eq:xi}, and $\uu={\uu}_k$ the fixed point of $\xi$ and $\xx=\xx_k=\pi(\uu)$. The language of $\xx$ inherits some symmetry from the language of $\uu=\xi^{\omega}(0)$. Once again, we use the notion of twin.

\begin{observation}
    Since the language of $\uu$ is closed under the letter permutation $\ell \leftrightarrow \ell'$ for each $\ell \in \{0,1,\dots, k-1\}$ and since $\pi(\ell') = 1-\pi(\ell)$ for each $\ell \in \{0,1,\dots, k-1\}$,  the language of $\xx$ is closed under the letter exchange  $0 \leftrightarrow 1$.  If $w'$ is obtained by this letter exchange from $w \in \mathcal{L}(\xx)$, we call $w'$ a {\em twin} of $w$. 
\end{observation}

First, we will examine for how long two words can have the same projection without matching letters in the same positions. An example of such words  $u,v\in \mathcal{L}(\uu)$ with  $\pi(u)=\pi(v)$ are  $u=u_1u_2u_3 = 201$ and $v=v_1v_2v_3 = 0'1'0'$. 

We will use the following observation. 

\begin{observation}\label{obs:uniquePreimage} 
    Let $w\in \mathcal{L}(\uu)$. If   $\pi(w) = 0101$,  then $w = 0101$.     
\end{observation}

\begin{proof} 
Denote $w = w_1w_2w_3w_4$. Using Observation~\ref{obs:basic} and the definition of $\pi$,  we have $w_3w_4=01$ or $w_4=0'$. Similarly, we find that either $w_2w_3=0'1'$ or $w_3=0$, and that either $w_1w_2=01$ or $w_2=0'$. Now, we use Observation~\ref{obs:basic} to rule out some combinations of those cases:
\begin{itemize}
    \item If $w_3w_4=00'$, then Item \ref{item:predek0} of Observation \ref{obs:basic} implies $w_2=(k-1)'$, a~contradiction.
    \item If $w_2w_3w_4=0'1'0'$, then it must be the image of $0'0'$ and must thus be preceded by the image of $(k-1)$ by the same Item, so $w_1=0'$, a~contradiction.
    \item If $w_2w_3w_4=0'01$, then the same Item forces again $w_1=(k-1)$, a~contradiction.
\end{itemize}
The only remaining possibility is $w=0101$, which concludes the proof.
\end{proof}

\begin{lemma}\label{lem:NoletterConicides} 
    Let $u =u_1u_2 \cdots u_n$ and  $v =v_1v_2 \cdots v_n$  belong to $\mathcal{L}(\uu)$.   Assume that $\pi(u) = \pi(v)$  and the letters  $u_\ell, v_\ell \in \mathcal{A}$ satisfy $u_\ell\neq v_\ell$  for every $\ell = 1,2, \ldots,n$. Then $n \leq k+1$. 

    If, moreover, $u0$ and $v0$ belong to $\mathcal{L}(\uu)$,  or $u0'$ and $v0'$ belong to $\mathcal{L}(\uu)$,  then $n\leq k-1$  and $\pi(u) = \pi(v)$  equals $1^{n}$ or $0^{n}$.
\end{lemma}

\begin{proof}
First assume that the last letters  $u_n$  and $v_n$ belong to $\{1,2,\ldots, k-1\}$. Say $u_n =i$ and  $v_n = j$, where $1\leq i < j\leq k-1$.  Item~\ref{item:j} of Observation \ref{obs:basic} gives that  $\pi(u) = \pi(v)$ is a suffix of $1^{i}$.  Clearly, $n=|\pi(u)|< k-1$.  
Analogously, if $u_n$  and $v_n$ belong to $\{1',2',\ldots, (k-1)'\}$,  we get that  $\pi(u) = \pi(v)$ is a suffix of $0^{i}$ for some $i<k-1$. 

\medskip

Let $u_n = i\in  \{1,2,\ldots, k-1\}$ and $v_n = 0'$. 

\begin{description}
\item {\bf Case} $\pi(u_n)= \pi(u_{n-1})=1$.\\
\noindent We consider $n\geq 3$ (if $n\leq 2$, there is no need to prove anything). By Observation \ref{obs:basic}, necessarily $i\geq 2$ and $v_{n-1} \in \{0',\ell\}$ for some $\ell\in \{1,2,\dots, k-1\}$, $\ell\neq i-1$. 
The same observation implies that 
\begin{itemize}
    \item $u_{n-2}u_{n-1}u_n = (i-2)(i-1)i$ \   and \  $\pi(u_{n-2}u_{n-1}u_n)$  is a suffix of $01^{i}$;

    \item $v_{n-2}v_{n-1}v_n =  (k-1)0'0'$, resp.  $(\ell-1)\ell 0'$  \ and  \ $\pi(v_{n-2}v_{n-1}v_n)$ is a suffix of $01^{k+1}$, resp.  \ $01^{\ell+1}$. 
\end{itemize}
Since $\pi(u)= \pi(v)$, we deduce that $\pi(u) = \pi(v) = 1^j$ for some $j \leq i$ and thus $n = |\pi(u)|=j \leq i\leq k-1$.

\item {\bf Case} $\pi(u_n)\neq  \pi(u_{n-1})=0$. \\    
\noindent Item~\ref{item:predek0} of Observation \ref{obs:basic}  implies that $v_{n-1}v_n$ equals $ 00'$ or $\ell' 0'$  for some $\ell\in \{1,2,\dots, k-1\}$.  Note that this situation cannot happen if  $v 0$ or $v0'$ is from  $\mathcal{L}(\uu)$, since by the same observation,  none of the factors \ $\ell' 0'0$, \  $ 0\,0' 0'$,  \ $\ell' 0'0'$, \   $ 0\,0' 0$  \  belong to the language of $\uu$.   Consequently, the second part of the lemma is already shown.  

To complete the proof, i.e., to show that $n\leq k+1$,  we assume that $n\geq 4$ (otherwise there is nothing to prove).  By Observation \ref{obs:uniquePreimage},   $\pi(u)$ cannot have a suffix $0101$.   Hence either $\pi(u_{n-1}) = \pi(u_{n-2})$ or $\pi(u_{n-2}) = \pi(u_{n-3})$.  
 
We apply the already proved statement to $\tilde{n} = n-1$ or to $\tilde{n} = n-2$ and we  deduce that   $\tilde{n}\leq k-1$. In particular,  $n\leq k+1$, as desired. 
\end{description}
\end{proof}

\begin{lemma}\label{lem:ShortCoincidenceToRight2}
    Let $w\neq \varepsilon$ be a bispecial factor of $\uu$. Assume that there exist $x,y \in \mathrm{Rext}(w)$  such that $x\neq y$ and $\pi(x)=\pi(y)$. Then $w$ or its twin equals $\xi^\ell(0)$, where $0\leq \ell \leq 2k-2$ and $\ell \neq k-1$. 

    Moreover, if $u , v \in \mathcal{A}^*$ are such that  $wxu, wyv \in \mathcal{L}(\uu)$  and $\pi(wxu)= \pi(wyv)$, then 
    \begin{itemize}
        \item $|u|= |v| \leq 1$, \  if \  $0\leq \ell \leq k-2$;
        \item $|u|= |v|  = 0$, \ if  \ $k\leq \ell \leq 2k-2$. 
    \end{itemize}
\end{lemma}

\begin{proof}
Recall that in Section \ref{subsec:triplets} we have obtained a complete description of the bispecial triplets in $\uu$ with a nonempty bispecial factor. Using this description, we may also find all possible right extensions of these factors. The following observations rule out many possibilities for $w$.
\begin{itemize}
    \item Considering Type 4, we note that the bispecial factor $w=\xi^{n+2k-1}(0)\xi^{n}(0)$ has two right extensions $x=0$ and $y=n+1$ for $n \in \{0,\dots, k-2\}$ and it has two right extensions $x=0$ and $y=0'$ for $n\geq k-1$, therefore $\pi(x)\not =\pi(y)$.
    \item Considering Types 1 to 3, we can see that the bispecial factor $w=\xi^n(0)$ for $n\geq 2k-1$ has two right extensions $x=0$ and $y=0'$, hence again $\pi(x)\not =\pi(y)$.
    \item Considering Type 1,  the right extensions of $w=\xi^{k-1}(0)$ are only $x=0$ and $y=0'$, thus once more $\pi(x)\not =\pi(y)$.
\end{itemize}
The first part of the statement is thus already proven, as the $x,y$ from the statement exist only when $w$ or its twin is $\xi^{\ell}(0)$ for $0\leq \ell\leq k-2$ or $k \leq \ell \leq 2k-2$.

To prove the second part of the statement, we examine the right extensions of these words. These extensions are listed in Table \ref{tab:bxu2}. Letters after a '$-$'character can be uniquely deduced from what stands before. We used the fact that the factors  $(k-1)'0'$  and $00$ in $\uu$ are always followed by $1'$ and $1$ respectively, compare with Items~\ref{item:(k-1)0} and~\ref{item:001} Observation \ref{obs:basic}. Note also that $s_R$ from Type 1 can be discarded as it is a prefix of the $s_R$ of Type 2. Similarly, $t_R$ from Type 1 can be discarded for $1\leq \ell\leq k-2$ as the other right extensions are all projected to the other letter, and for $k\leq \ell\leq 2k-2$ as it is equal to $s_R$ from Type 4.
Examining Table \ref{tab:bxu2}, we see that the right extensions of these $w$ cannot exceed the stated lengths while keeping equal projections, as expected.  
\end{proof}

\begin{table}[ht]
    \centering
    \setlength{\tabcolsep}{3pt}
    \renewcommand{\arraystretch}{1.2}
\begin{tabular}{|c|l|l|l|l|}
\hline
$w$ &  $ w\,t_R\!\!-\!\! \tilde{u}$ (Type 2) &$w\,s_R$ (Type 2) &$ w\, s_R$ (Type 3) &$w\, s_R\!\!-\!\!\tilde{v}$ (Type 4)\\
\hline \hline
$\varepsilon$ &$w\,(k\!\!-\!\! 1 )'\,0'\!\!-\!\!1'$ &$w\, 01 $ & $w\, 00' $&$w\, 00\!\!-\!\!1 $  \\
\hline 
$0$  &$w\,0' \!\!-\!\! 1'2'$ &  $w\,12$ & $w\,10'$ &  $w\,10\!\!-\!\! 1$\\
\hline 
$\xi(0)$  &$w\,0' \!\!-\!\! 1'2'$ &  $w\,23$ & $w\,20'$ &  $w\,20\!\!-\!\! 1$\\
\hline 
$\vdots$&$ \vdots$ & $\vdots$ & $\vdots$ & $\vdots$ \\
\hline 
$\xi^{k-2}(0)$&$w\,0' \!\!-\!\! 1'2'$ &  $w\,(k\!\!-\!\!1)0'$ & $w(k\!\!-\!\!1)0'$ &  $w(k\!\!-\!\!1)0\!\!-\!\!1 $\\
\hline 
$\xi^k(0)$ &  $ w\,1'2'$ & \phantom{$w\, 0'$}&  \phantom{$w\, 0'$}& $w\,0\!\!-\!\!1$\\
\hline 
$\xi^{k+1}(0)$ &  $ w\,2'3'$ &  \phantom{$w\, 0'$}&  \phantom{$w\, 0'$}& $w\,0\!\!-\!\!1$\\
\hline 
$\vdots$&$ \vdots$ &  & & $\vdots$ \\
\hline
$\xi^{2k-2}(0)$ & $ w\,(k\!\! - \!\!1)'0$&\phantom{$w\, 0'$}&  \phantom{$w\, 0'$}& $w\,0\!\!-\!\!1$\\
\hline 
\end{tabular}\caption{Right extensions  of bispecial factors in the case when the first letters of $s_R$  and $t_R$ have the same projection by $\pi$.}\label{tab:bxu2}
\end{table}

\begin{lemma}\label{lem:CoincidenceToLeft2}
    Let $\ell \in \N$ satisfy $2\leq \ell \leq k-2$ or $k+2\leq \ell \leq 2k-2$ and $w = \xi^\ell(0)$.
    Assume that 
    \begin{enumerate}
    \item $x,y \in \mathrm{Rext}(w)$ are such that $x\neq y$ and $\pi(x)=\pi(y)$  and   
    \item $u, v \in \mathcal{A}^*$ are such that  $uwx, vwy \in \mathcal{L}(\uu)$  and  $\pi(u) = \pi(v)$ is a prefix of $1^{k-1}$ or $0^{k-1}$. 
    \end{enumerate}
    Then $|u| = |v| \leq k-3$.     
\end{lemma}

\begin{proof}
For  bispecial factors $w=\xi^{\ell}(0)$ with  $\ell$ specified in the assumptions of the lemma, Table~\ref{tab:ubx2} shows all possible right extensions $x,y$ which satisfy Item 1. 

If $2\leq \ell \leq k-2 $, the table shows that the letter $(\ell-1)$ is the  left extension of the factor $w0'$. By Observation \ref{obs:basic}, the factor $w0'$  is preceded by  $01\cdots (\ell- 2)(\ell- 1)$. As $\pi\bigl(01\cdots (\ell- 2)(\ell- 1)\bigr) = 01^{\ell-1}$, the factor $u$ satisfying Item 2 of the lemma is of length at most  $\ell-1 \leq k-3$. 

If $k+2\leq \ell \leq 2k-2 $, Table \ref{tab:ubx2} says  that the factor $w0$ has the left extension $j':= (\ell\!\!-\!\! k\!\!-\!\! 1)'$, and therefore $w0$ is preceded by the factor $0'1'\cdots (j\!\!-\!\! 1)'j'$. As $\pi\bigl(0'1'\cdots (j\!\!-\!\! 1)'j'\bigr) = 10^{j}$, the factor $v$ satisfying Item 2 is
of length at most $j=\ell-k-1 \leq  k-3$. 
\end{proof}

\begin{table}[ht]
    \centering
    \setlength{\tabcolsep}{3pt}
    \renewcommand{\arraystretch}{1.2}
\begin{tabular}{|c|c|c|}
\hline
$w$ &  ~~~$t_L\, w\, t_R\!\!-\!\! \tilde{u}$ &~~~$s_L\,w\,y $\\
\hline \hline
$\varepsilon$ &   $~~(k\!\! - \!\!2)'\, w\, (k\!\!-\!\!1)'0'$ & $~~~~~*\ w\,0~~~$\\
\hline 
$0$ &$(k\!\!-\!\! 1 )'\, w\, 0'\!\!-\!\! 1'~~~$ & $~~~~~* \  w\, 1~~$  \\
\hline 
$\xi(0)$ & $~~~~~~0\,\ w\, 0'\!\!-\!\! 1'~~~$  &$~~~~~* \ w \,  2~~$\\
\hline 
$\xi^2(0)$  & $~~~~~~1\,\ w\, 0'\!\!-\!\! 1'~~~$ &$~~~~~* \ w\,    3~~$\\
\hline 
$\vdots$&$ \vdots$ & $\vdots$  \\
\hline 
$\xi^{k-2}(0)$&  $(k\!\!-\!\! 3)\  w\, 0'1'~~~$& ~~$~~~~~~*\ w\,   (k\!\!-\!\! 1)$\\
\hline 
$\xi^k(0)$ & $(k\!\! - \!\! 2 )\ w\,1'2'~~~$ &  $(k\!\! -\!\! 1)\ w\,0~~~$\\
\hline 
$\xi^{k+1}(0)$ &  $(k\!\! - \!\! 1 )\ w\, 2'3'~~~$ & $~~~~~0'\ w\, 0~~~$\\
\hline 
 
$\xi^{k+2}(0)$ &  $~~~~0'\ w\, 3'4'~~~$ & $~~~~~1' \ w\,0~~~$\\
\hline 
$\vdots$&$ \vdots$ & $\vdots$\\

\hline
$\xi^{2k-2}(0)$ &   $(k\!\! - \!\!4)'\, w\, (k\!\!-\!\!1)'0$ & $(k\!\!-\!\! 3)'\, w\,0~~~$\\
\hline 
\end{tabular}\caption{Left and right extensions  of bispecial factors $\xi^{\ell}(0)$ for $0\leq\ell\leq 2k-2, \ \ell \not =k-1$, where $y$ stands for the first letter of $s_R$,  $*$ is any  letter from  $ \mathcal{A}$ different from $t_L$. Notation   $t_L\, w\, t_R \!-\! \tilde{u}$  means that the factor  $t_L\, w\, t_R$ is in $\uu$ always followed by $\tilde{u}$.}\label{tab:ubx2}
\end{table}

\section{Proof of Conjecture \ref{conj:complexity}}\label{sec:proof_complexity}

We have prepared everything to reveal a close relation between left special factors in $\uu$ and $\xx$, which is crucial for the proof of Conjecture~\ref{conj:complexity}.
    
\begin{proposition}\label{pro:LSinX}
    Let $w$ be a left special factor of $\xx$  such that  $|w|$ is at least $ |\xi^{2k-2}(0)|+k+2$. 
    Then there exist $i\in \N, 1\leq i\leq k$, and a left special factor $f$ of $\uu$ such that $w$ or its twin equals $w=0^{i-1}\pi(f)$.   
\end{proposition}

\begin{proof} Denote  $n = |w|$. As $0w$ and $1w$ belong to $\mathcal{L}(\xx)$, there exist factors $u=u_1u_2\cdots u_n$  and $v =v_1v_2 \cdots v_n$ and letters $a,b$ such that $au, bv \in \mathcal{L}(\uu)$ 
and  $0w = \pi(au)$ and $1w=\pi(bv)$. 
We denote 
$$
i =\min\{\ell : u_\ell = v_{\ell}\}\qquad \text{and}\qquad 
j = \max\{\ell: u_iu_{i+1}\cdots u_{\ell} = v_iv_{i+1}\cdots v_{\ell}\}. 
$$
By Lemma~\ref{lem:NoletterConicides}, $i\leq  k+2\leq n$.  Thus the letter $u_{i} = v_{i}$ has two left extensions. In particular, by Observation~\ref{obs:basic}, $u_i \in \{0,0'\}$. Again, by Lemma~\ref{lem:NoletterConicides}, second part, $i\leq k$  and $\pi(u_1u_2\cdots u_{i-1}) = 1^{i-1}$ or $0^{i-1}$.  The definition of $i$ implies that the factor $f =u_iu_{i+1}\cdots u_j = v_iv_{i+1}\cdots v_j$ is a left special factor in $\mathcal{L}(\uu)$. To prove the proposition we just need to prove that 
$j=n$.
  
Assume the contrary, i.e., that $j< n$. 
Then $f$ is also a right special factor with two right extensions $u_{j+1}\neq v_{j+1}$ and $\pi(u_{j+1}) = \pi(v_{j+1})$. By  Lemma~\ref{lem:ShortCoincidenceToRight2},  the length  $j-i+1$ of $f$ is at most $|\xi^{2k-2}(0)|$. This implies 
$$
n\geq |\xi^{2k-2}(0)| +k+2 \geq j-i+1+k+2 \geq j+3.
$$

By the same lemma, $\pi(u_{j+2})\neq \pi(v_{j+2})$ or $\pi(u_{j+3})\neq \pi(v_{j+3})$, which contradicts the fact that $\pi(u)=\pi(v)$.  
\end{proof}

\begin{proposition}\label{pro:BSinX} 
    Let $f \in \mathcal{L}(\uu)$ be a left special factor of $\uu$ such that $0$ is the first letter of $f$. 
    \begin{enumerate}
    \item If $f$ is a prefix of $\uu$, then $0^{i-1}\pi(f)$ and $1^{i-1}\pi(f)$ are left special factors of $\xx$ for every $i\in \{1,2,\dots, k\}$.
    
    \item If $f$ is not a prefix of $\uu$,  then there exist a unique  $c\in \{0,1\}$ and a unique $i \in \{1,2,\ldots,k\}$ such that $w=c^{i-1}\pi(f)$ is a left special factor  in $\xx$.
    \end{enumerate}
\end{proposition}

\begin{proof} 
Let us recall that by Observation~\ref{obs:basic} the following hold for every $j \in \{1,2,\ldots, k-1\}$.

\begin{itemize}
    \item The factor $j0$ occurs  in $\uu$ only as a suffix  of the factor $01\cdots (j-1)j0$.
    \item The factor $j'0$ occurs  in $\uu$ only as a suffix   of the factor  $0'1'\cdots (j-1)'j'0$.
    \item The factor $00$ occurs in $\uu$ only as  a suffix  of the factor  $0'1'\cdots (k-1)'00$.
    \item The factor $0'0$ occurs in $\uu$ only as  a suffix  of the factor  $01\cdots (k-1)0'0$.
\end{itemize}


\noindent Let  $a \in \mathcal{A}$ be a left extension of $f$. Since $f$ starts with $0$,  the factor  $\pi(af)$ occurs  in $\xx$  only as a suffix of 
\begin{itemize}
    \item $01^j\pi(f)$, \  if $a=j\in \{1,\ldots, k-1\}$,
    \item $10^j\pi(f)$, \  if $a=j'\in \{1',\ldots, (k-1)'\}$,
    \item $10^k\pi(f)$, \  if $a=0$,
    \item $01^k\pi(f)$, \  if $a=0'$. 
\end{itemize}

\noindent Let us stress that to  distinct letters that could be left extensions of $f$, we assigned  distinct factors in the previous list. 
\medskip

\noindent Let $w \in \mathcal{L}(\uu)$ be the shortest bispecial factor having $f$ as its prefix. 
\begin{enumerate}
    \item If $f$ is a prefix of  $\uu$, the bispecial factor $w$ equals $\xi^n(0)$  for some $n \in \N$ and thus $\mathrm{Lext}(w) = \mathrm{Lext}(f) = \mathcal{A}$ due to Proposition \ref{prop:aboutBSfactors}. Using the previous list, we see that for every $\ell\in \{1,2\ldots, k\}$, 
    \[01^{\ell-1} \pi(f),\text{ }
    11^{\ell-1} \pi(f),\text{ }
    10^{\ell-1} \pi(f),\text{ and }
    00^{\ell-1} \pi(f)\]
    all belong to $\mathcal{L}(\xx)$. Therefore,  $1^{\ell-1}\pi(f)$  and $0^{\ell-1}\pi(f)$  are  left special factors in $\xx$. 

    \item If $f$ is not a prefix of  $\uu$, the bispecial factor $w$ equals $\xi^{n+2k-1}(0)\xi^n(0)$  for some $n \in \N$ and thus $\#\mathrm{Lext}(w) = \#\mathrm{Lext}(f) = 2$. Then $\pi(f)$ preceded by the longest common suffix of the factors corresponding to the two letters from $\mathrm{Lext}(f)$ in the above list is a left special factor in $\xx$. 
\end{enumerate}
\end{proof}

The following theorem confirms the validity of  Conjecture~\ref{conj:complexity}.

\begin{theorem}\label{thm:complexity}
    For sufficiently large $N$, the first difference of factor complexity of the sequence $\xx$ satisfies 
    $$
    \Delta {\mathcal C}(N)\in \{4k-2, 4k\}\,.
    $$
\end{theorem}

\begin{proof} 
By Formula \eqref{eq:ComplexityDifference}, we have  to determine the number of left special factors in $\xx$  of length $N$. 

Since $\uu=\xi^\omega(0)$, the bispecial factor $\xi^n(0)$ is a prefix of $\uu$ for every $n \in \N$ and thus any prefix of $\uu$ is a left special factor in $\uu$.  By Proposition \ref{pro:BSinX}, for every $i\in \{1,2,\ldots, k\}$  and every prefix $f$ of $\uu$ of length $N-i+1$, the words $0^{i-1}\pi(f)$ and $1^{i-1}\pi(f)$ are left special factors of length $N$ in  $\xx$. As the language of $\xx$ is closed under exchange $0\leftrightarrow 1$, the same fact is true for the twins of the above mentioned left special factors. It is easy to see that for  $N>k$, all these left special factors differ. Thus, for every $N$ we have at least $4k-2$ left special factors in $\xx$.

Let us consider a left special factor $f$ of  $\uu$  which is not a prefix of $\uu$ and such that $0$ is the first letter of $f$.
By Proposition~\ref{prop:aboutBSfactors}, the shortest bispecial factor of $\uu$ with the prefix $f$  has the form $v=\xi^{n+2k-1}(0)\xi^n(0)$. We will use the following  properties of $v$.  

\begin{itemize}
    \item The longest common prefix of the bispecial factor $v=\xi^{n+2k-1}(0)\xi^{n}(0)$  and $\uu$ is  $w=\xi^{n+2k-1}(0)$.

Indeed, on the one hand $w0$ is a prefix of $v$ as $0$ is a prefix of $\xi^n(0)$, and on the other hand we show that $w0'$ is a prefix of $\xi^{n+2k}(0)$ which is a prefix  of $\uu$. We have
$$
\xi^{n+2k}(0)= \xi^{n+2k-1}(01) = w\xi^{n+2k-1}(1)\,,  
$$ 
where by Item~\ref{item:1} of Observation~\ref{obs:basic}, $0'$ is a prefix of $\xi^{n+2k-1}(1)$. 
     
    \item  $|v|+k-1<|\xi^{n+2k}(0)|$ for every $n\in \mathbb N$.
      
Indeed,        
$$
\begin{array}{rcl}
    |v|+k-1&=&|\xi^{n+2k-1}(0)\xi^{n}(0)|+k-1\\
    &\leq& |\xi^{n+2k-1}(0)\xi^{n}(01\cdots (k-1))|\\
    &<& |\xi^{n+2k-1}(0)\xi^{n+k}(0)|=|\xi^{n+2k-1}(0)\xi^{n+k}(0')|\\
    &=& |\xi^{n+2k-1}(0)\xi^{n+2k-1}(1)|=|\xi^{n+2k}(0)|\,.
\end{array}
$$ 
\end{itemize}

We use the  first property of $v=\xi^{n+2k-1}(0)\xi^n(0)$ mentioned above and  the constant $c$ from Item 2 of Proposition \ref{pro:BSinX}.  They  imply  that for every $N$ such that   $ |\xi^{n+2k-1}(0)| + c \leq N \leq  |v| + c-1$, we have   a new left special factor of length $N$.  Since its twin is left special, too, we have to add 2 to the number of left special factors of length $N$. 

For every $n \in \N$, we have to consider  
left special factors  in $\xx$  arising from the bispecial factor $\xi^{n+2k-1}(0)\xi^n(0)$. The second property of $v$ implies that the lengths  of  left special factors arising from  $\xi^{n+2k-1}(0)\xi^n(0)$ are strictly smaller than the length of those arising from  $\xi^{n+2k}(0)\xi^{n+1}(0)$. 

\medskip

Proposition \ref{pro:LSinX} confirms that we have already accounted for all possible left special factors of length $N$.  

\end{proof}

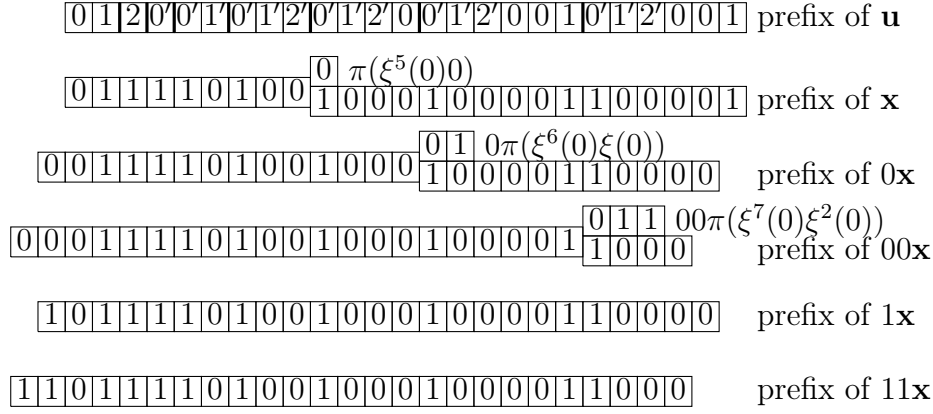
\begin{figure}
    \centering
\begin{tikzpicture}[xscale=0.36,yscale=0.395]
    \foreach\x [count=\xi] in {0,1,2,0',0',1',0',1',2',0',1',2',0,0',1',2',0,0,1,0',1',2',0,0,1}{
            \draw (\xi,0)--(\xi,1)--(\xi+1,1)--(\xi+1,0)--(\xi,0);
            \node at (\xi+0.5,0.5){$\strut\x$};};
    \draw[very thick] (3,1)--(3,0);
    \draw[very thick] (4,1)--(4,0);
    \draw[very thick] (5,1)--(5,0);
    \draw[very thick] (7,1)--(7,0);
    \draw[very thick] (10,1)--(10,0);
    \draw[very thick] (14,1)--(14,0);
    \draw[very thick] (20,1)--(20,0);
    \foreach\image[count=\rownumber] in 
    {
        {0,1,1,1,1,0,1,0,0},
        {0,1,1,1,1,0,1,0,0,1,0,0,0},
        {0,1,1,1,1,0,1,0,0,1,0,0,0,1,0,0,0,0,1},
        {0,1,1,1,1,0,1,0,0,1,0,0,0,1,0,0,0,0,1,1,0,0,0,0},
        {0,1,1,1,1,0,1,0,0,1,0,0,0,1,0,0,0,0,1,1,0,0,0}}
    {
        \foreach\x [count=\xi] in \image{
            \draw (\xi,-2.5*\rownumber)--(\xi,-2.5*\rownumber+1)--(\xi+1,-2.5*\rownumber+1)--(\xi+1,-2.5*\rownumber)--(\xi,-2.5*\rownumber);
            \node at (\xi+0.5,-2.5*\rownumber+0.5){$\strut\x$};};
    }
    \draw(0,-5)--(0,-4)--(1,-4)--(1,-5)--(0,-5);
    \node at (0.5,-4.5){$\strut 0$};
    \draw(0,-7.5)--(1,-7.5)--(1,-6.5)--(0,-6.5)--(0,-7.5);
    \node at (0.5,-7){$\strut 0$};
    \draw(-1,-7.5)--(0,-7.5)--(0,-6.5)--(-1,-6.5)--(-1,-7.5);
    \node at (-0.5,-7){$\strut 0$};
    \draw(0,-10)--(1,-10)--(1,-9)--(0,-9)--(0,-10);
    \node at (0.5,-9.5){$\strut 1$};
    \draw(0,-12.5)--(1,-12.5)--(1,-11.5)--(0,-11.5)--(0,-12.5);
    \node at (0.5,-12){$\strut 1$};
    \draw(-1,-12.5)--(0,-12.5)--(0,-11.5)--(-1,-11.5)--(-1,-12.5);
    \node at (-0.5,-12){$\strut 1$};
    \begin{scope}[xshift=9cm,yshift=-0.8cm]
    \foreach\image[count=\rownumber] in 
    {
        {0},
        {1,0,0,0,1,0,0,0,0,1,1,0,0,0,0,1}}
    {
        \foreach\x [count=\xi] in \image{
            \draw(\xi,-\rownumber)--(\xi,-\rownumber+1)--(\xi+1,-\rownumber+1)--(\xi+1,-\rownumber)--(\xi,-\rownumber);
            \node at (\xi+0.5,-\rownumber+0.5){$\strut\x$};};
    }
    \end{scope}
    \begin{scope}[xshift=13cm,yshift=-3.3cm]
    \foreach\image[count=\rownumber] in 
    {
        {0,1},
        {1,0,0,0,0,1,1,0,0,0,0}}
    {
        \foreach\x [count=\xi] in \image{
            \draw(\xi,-\rownumber)--(\xi,-\rownumber+1)--(\xi+1,-\rownumber+1)--(\xi+1,-\rownumber)--(\xi,-\rownumber);
            \node at (\xi+0.5,-\rownumber+0.5){$\strut\x$};};
    }
    \end{scope}
    \begin{scope}[xshift=19cm,yshift=-5.8cm]
    \foreach\image[count=\rownumber] in 
    {
        {0,1,1},
        {1,0,0,0}}
    {
        \foreach\x [count=\xi] in \image{
            \draw(\xi,-\rownumber)--(\xi,-\rownumber+1)--(\xi+1,-\rownumber+1)--(\xi+1,-\rownumber)--(\xi,-\rownumber);
            \node at (\xi+0.5,-\rownumber+0.5){$\strut\x$};};
    }
    \end{scope}
    \node [anchor=west] at (26,0.5) {prefix of $\uu$};
    \node [anchor=west] at (11,-1.3) {$\pi(\xi^5(0)0)$};
    \node [anchor=west] at (26,-2.3) {prefix of $\xx$};
    \node [anchor=west] at (16,-3.8) {$0\pi(\xi^6(0)\xi(0))$};
    \node [anchor=west] at (26,-4.8) {prefix of $0\xx$};
    \node [anchor=west] at (23,-6.3) {$00\pi(\xi^7(0)\xi^2(0))$};
    \node [anchor=west] at (26,-7.3) {prefix of $00\xx$};
    \node [anchor=west] at (26,-9.5) {prefix of $1\xx$};
    \node [anchor=west] at (26,-12) {prefix of $11\xx$};
\end{tikzpicture}
    \caption{All left special factors of $\xx_k$, where $k=3$, of length $N$ with $10\leq N \leq 25$, are prefixes of the depicted factors or their twins.}
    \label{fig:leftspecial-branches}
\end{figure}

The structure of left special factors is illustrated in Figure~\ref{fig:leftspecial-branches}. In the proofs, we limit our consideration to $k\geq 4$ because Conjecture~\ref{conj:complexity} has been already proven for $k\leq 3$. Nevertheless, due to space constraints,  Figure 1  exhibits the situation for $k=3$.

\section{Proof of Conjecture \ref{conj:repetitions}}\label{sec:proof_repetitions}

Using the fact that $\uu$ contains no overlaps, see Proposition~\ref{pro:NoOverlap}, and our knowledge of the left and right extensions of bispecial factors in $\uu$, see~ Lemma~\ref{lem:ShortCoincidenceToRight2}, we are ready to prove Conjecture~\ref{conj:repetitions}.

\begin{theorem}\label{thm:repetitions}
    The sequence $\xx$ has critical exponent $k + 1$, which is attained by the words $0^{k+1}$ and ${1}^{k+1}$. It contains no factor of length $2n + k$ and period $n$, and therefore has asymptotic critical exponent equal to $2$. 
\end{theorem}

\begin{proof}  
By Observation \ref{obs:basic}, the words $1^{k+1}, 0^{k+1}$ are in $\mathcal{L}(\xx)$,  but $1^{k+2}, 0^{k+2}$ are not.   Therefore, we consider $n\geq 2$. 

Assume the contrary to the statement of the theorem, i.e., that there exists in $\xx$ an overlap of size $k$. Equivalently, denoting  $m=n+k$, there exists a  factor  of $\uu$ of length $2m-k$,  say, 
$$
w=w_1w_2 \cdots w_{m-k}w_{m-k+1} \cdots w_mw_{m+1} \cdots w_{2m-k} \in \mathcal{L}(\uu)
$$ 
such that  $m \geq k+2$ and  the projection by $\pi$ of the prefix of $w$ of length $m$ coincides with the projection by $\pi$ of the suffix of $w$ of length $m $, formally, 
\begin{equation}\label{eq:overlapInProjection}
    \pi(w_1w_2 \cdots w_{m-k}w_{m-k+1} \cdots w_m) = \pi(w_{m-k+1} \cdots w_mw_{m+1} \cdots w_{2m-k})\,.
\end{equation}
Let us denote 
$$
i =\min\{\ell : w_\ell = w_{m-k+\ell}\}.
$$
Lemma  \ref{lem:NoletterConicides} implies $i \leq k+2$. If $i$ were at least $2$, then   the letter $w_i = w_{m-k+i}$ would have  two distinct letters in its left extension, namely $w_{i-1}$ and $w_{m-k+i-1}$. Hence, $w_i$ is either $0$ or $0'$. Put $u = w_1w_{2}\cdots w_{i-1}$ and $v=w_{m-k+1}w_{m-k+2}\cdots w_{m-k+i-1}$. Clearly, $\pi(u)=\pi(v)$.  Since $u0', v0' \in \mathcal{L}(\uu)$ or $u0, v0 \in \mathcal{L}(\uu)$, Lemma \ref{lem:NoletterConicides} implies that $i-1\leq k-1$. We have shown 
\begin{equation}\label{eq:iBig}
    \text{either} \quad i=1\quad \text{or}\quad i\leq k.
\end{equation}

Denote 
$$
j = \max\{\ell: w_iw_{i+1}\cdots w_{\ell} = w_{m-k+i}w_{m-k+i+1}\cdots w_{m-k+\ell}\}\,. 
$$
Obviously, $i\leq j  \leq m$. 
If $j$ were at least $m-k+i$, then $w_iw_{i+1}\cdots w_{m-k+j}$ would be an overlap in $\uu$, which contradicts  Proposition \ref{pro:NoOverlap}. From now on, we assume that 
\begin{equation}\label{eq:onlyParametersFor discussion}
    j \leq m-k+i-1. 
\end{equation}
If $j$ were equal to $m$, then the previous inequality would force  $i\geq k+1$ -- this contradicts  \eqref{eq:iBig}. 

\medskip

Necessarily,  $j\leq m-1$. 
The definition of $j$ forces $w_{j+1}\neq w_{m-k+j+1}$. In particular, 
$$
w_iw_{i+1}\cdots w_{j} = w_{m-k+i} w_{m-k+i+1} \cdots w_{m-k+j} \quad \text{is a right special factor in $\uu$},
$$
with two extensions $w_{j+1}\neq w_{m-k+j+1}$ having the same projection by $\pi$. This means that $f:=w_iw_{i+1}\cdots w_{j}$ is a suffix of a bispecial factor  $\xi^\ell(0)$ or its twin considered in Lemma~\ref{lem:ShortCoincidenceToRight2}. 

Due to \eqref{eq:overlapInProjection}, we also have 
$$
\pi(w_iw_{i+1} \cdots  w_m) = \pi(w_{m-k+i}w_{m-k+i+1} \cdots  w_{2m-k}). 
$$
Using Lemma~\ref{lem:ShortCoincidenceToRight2} with $w=w_iw_{i+1}\cdots w_{j}=w_{m-k+i} w_{m-k+i+1} \cdots w_{m-k+j}$ and $xu$ prefix of $w_{j+1}\cdots w_m$ and $yv$ prefix of $ w_{m-k+j+1}\cdots w_{2m-k}$, since $|u|=|v|\leq 1$, we have $j+2\geq m$, which leads to $j \geq m-2$.
\medskip

If $i$ were at most $k-2$, then $j \geq m-2 = m-k +k-2 \geq m-k +i$, which again  contradicts \eqref{eq:onlyParametersFor discussion}. 
 From now on we may restrict our consideration to the values $i, j$ such that 
\begin{equation}\label{eq:narrowParameters}
    k-1\leq  i \leq k \qquad \text{and} \qquad m-2\leq j \leq m-1 .
\end{equation}
Since $k\geq 4$, we have $i\geq 3$.  The  definition of $i$ implies $w_{i-1} \neq w_{m-k+i-1}$. Therefore, the factor 
\begin{equation}\label{eq:f}
    f=w_i\cdots w_j = w_{m-k+i}\cdots w_{m-k+j} \quad \text{is bispecial}. 
\end{equation}  
More precisely,   $f$ or its twin equals  $\xi^{\ell}(0)$ for some $\ell \in \{0,1,\ldots, 2k-2\} \setminus \{k-1\}$. 

As $\pi(w_1 \cdots w_jw_{j+1}) = \pi(w_{m-k+1} \cdots w_{m-k+j}w_{m-k+j+1})$,  one can derive  $\pi(ufx) =  \pi(vfy)$ denoting   
$$
u =w_1 \cdots w_{i-1}, \  x = w_{j+1} \ \text{ and}\ \  v = w_{m-k+1} \cdots w_{m-k+i-1},  \ y = w_{m-k+j+1}.
$$
Note that $\pi(u)$ is $0^{i-1}$ or $1^{i-1}$ from Lemma \ref{lem:NoletterConicides}. Under the assumption $\ell \not \in\{0,1, k, k+1\}$,  Lemma~\ref{lem:CoincidenceToLeft2} forces $|u| = |v| =i-1 \leq k-3$. This contradicts \eqref{eq:narrowParameters}. Consequently, it remains to treat $\ell \in\{0,1, k, k+1\}$.

\medskip

Let us recall that  $f=w_i \cdots w_j$  and the index $m-k+i$ is also an  occurrence of $f$ in $w$.   By \eqref{eq:onlyParametersFor discussion} and \eqref{eq:narrowParameters}
$$
j+1\leq m-k+i = \underbrace {m-2-j}_{\leq 0}+j+\underbrace{i-k}_{\leq 0}+2 \leq j+2, 
$$
i.e.,  the first occurrence of $f$ ends at position $j$  and simultaneously $f$  occurs again at position $j+1$ or $j+2$. 
Note that $z_1:=w_{j+1}$ and $z_2: = w_{m-k+j+1}$ belong to the right extensions of $f$, $z_1\neq z_2$ and $\pi(z_1)=\pi(z_2)$. 
Hence, 
\begin{description}
\item[either (A):] \quad  the index $i$ is an occurrence of  $ffz_2$ and $z_1$ is a prefix of $f$;
\item[or (B):]  \quad the index $i$ is an occurrence of $fz_1fz_2$.    
\end{description}
 We look at Table \ref{tab:ubx2} to find $\{z_1,z_2\}$, to determine what precedes  $fz$ if $z \in \{z_1,z_2\}$ and to determine what follows $fz$. 
\begin{description}
\item[Case $f= 0$:] \quad  $\{z_1,z_2\}=\{0',1\}$. As none of $z_1$ and $z_2$ is a prefix of $f$, the situation $(A)$ is excluded.
 By the table, the factor  $f0' = 00'$ is preceded by $(k-1)'$ and  followed by $1'$. Both possibilities for (B), namely $fz_1fz_2 = f0'f1$ and $fz_1fz_2 = f1f0'$ lead to a contradiction. 

\item[Case $f = \xi(0)=01$:] \quad $\{z_1,z_2\}=\{0',2\}$ and again $(A)$ is excluded as $0$ is a~prefix of $f$.   By Table~\ref{tab:ubx2}, the factor  $f0' = \xi(0)0'$ is preceded by $0$ and  followed by $1'$.  Both possibilities  $fz_1fz_2 = f2f0'$ and $fz_1fz_2 = f0'f2$ lead to a contradiction. 
\item[Case $f = \xi^{k}(0)$:] \quad  $\{z_1,z_2\}=\{0,1'\}$. By Table \ref{tab:ubx2}, the factor $f1' = \xi^{k}(0)1'$ is followed by $2'$.  Hence,  $z_1\neq 1'$, i.e., $z_2=1'$. By the table, the factor  $fz_2 = \xi^{k}(0)1'$ is preceded by $(k-2)$. This yields a~contradiction,  as the factor $f =  \xi^{k}(0)$ has a suffix $0'$ and $z_1 = 0$.

\item[Case $f = \xi^{k+1}(0)$:] \quad  $\{z_1,z_2\}=\{0,2'\}$.  By Table \ref{tab:ubx2}, the factor $f2' = \xi^{k+1}(0)2'$ is followed by $3'$.  Hence,  $z_1\neq 2'$, i.e., $z_2=2'$. But  the factor  $fz_2 = \xi^{k+1}(0)2'$ is preceded by $(k-1)$. This yields a contradiction,  as the factor $f =  \xi^{k+1}(0)$ has a suffix $1'$ and $z_1 = 0$.
\end{description}

We have proved that the assumption that there exists an overlap of size $k$ in $\xx$ leads to a contradiction. Therefore,  Conjecture 2 is proven.  
\end{proof}

\section{Comments and open problem}

\begin{itemize}
    \item It is well known that the factor complexity $\mathcal{C}$ of a coding of a fixed point of a primitive morphism is sublinear, i.e., there exists a constant $c$ such that $\mathcal{C}(n) \leq c\, n$ for every $n \in \N$. The validity of Conjecture 1 shows that the constant  $c$ cannot be bounded even for binary words. 

    \item  Binary infinite words can have an arbitrarily large difference between the critical exponent and the asymptotic critical exponent. Sturmian sequences are examples of such sequences with language closed under reversal. Their (asymptotic) critical exponent may be computed using the continued fraction expansion coefficients of their slope~\cite{DaLe}: while the critical exponent is larger than any coefficient, sequences sharing the same suffix for the continued fraction of their slope have the same asymptotic critical exponent. The validity of Conjecture 2 implies that the Thue-Morse-like sequences serve as another example, this time of sequences with language closed under letter exchange. 

    \item Let us fix $k \in \N$ and denote $x_0x_1x_2\ldots $ the Thue-Morse-like sequence  $\xx_k = \pi(\uu_k)$ studied in this paper.  In the case $k =3$,  Theorem 43 in  \cite{Sha25} states  that every $n\in \N$ can be written as a sum of two integers from the set  $J_0 = \{i : x_i =0 \}$ and also as a sum of two integers from the set $J_1 = \{i : x_i =1 \}$. Our computer experiments suggest that this statement is valid for every sequence $\xx_k$ with $k \geq 2$, although it is known to be false for the Thue-Morse sequence $\xx_1$ \cite{RSS20}. 
\end{itemize}




\bibliographystyle{elsarticle-num} 
\bibliography{biblio_ShallitNC.bib}



\end{document}